\documentclass[12pt,a4paper,titlepage]{article}

\usepackage{amssymb}

\usepackage{amsmath} \usepackage{graphicx} \usepackage{wrapfig} 

  \def\N{{\mathbb N}}     \def\R{{\mathbb R}} \def\pf{\emph{Proof. }} \def\qed{$\hfill\blacksquare$}  \def\cone{{\rm cone\,}}    \def\ocone{\overline{\rm cone}\,} \def\*ocone{^{w^*}\overline{\rm cone}\,}                                          \def\spane{{\rm span}\,}  \def\la{\langle} \def\ra{\rangle}  \def\alp{\alpha}    \def\del{\delta} \def\cP{\mathcal P}   \def\cB{\mathcal{B}} \def\cN{\mathcal{N}} \def\cR{\mathcal{R}} \def\cD{\mathcal{D}} \def\cone{{\rm cone}\,}    \def\be{\begin{equation}} \def\ee{\end{equation}}   \newtheorem{thm}{Theorem}[section] \newtheorem{lem}[thm]{Lemma}  \newtheorem{cor}[thm]{Corollary}  \newtheorem{defn}[thm]{Definition} \newtheorem{ex}[thm]{Example} \newtheorem{rem}[thm]{Remark}  \newtheorem{fact}[thm]{Fact}  \numberwithin{equation}{section} 

\begin{document} 

\title{Some Applications of the Hahn-Banach Separation Theorem \thanks{ We dedicate this paper to our friend and colleague Professor Yair Censor on the occasion of his 75th birthday.}} \vspace{.1in} \author{Frank Deutsch ,  Hein Hundal,  Ludmil Zikatanov \\} 

\maketitle                   \begin{abstract} We show that a single special separation theorem (namely, a consequence of the geometric form of the Hahn-Banach theorem) can be used to prove  Farkas type theorems, existence theorems for numerical quadrature with positive coefficients, and detailed characterizations of best approximations from certain important cones in Hilbert space. \vspace{12pt} 

 \end{abstract} 

\section{Introduction}\label{S:intro} 

We show that a single separation theorem---the geometric form of the Hahn-Banach theorem---has a variety of different applications.  

In section \ref{S: key thm} we state this general separation theorem (Theorem \ref{big thm}), but note that only a special consequence of it is needed for our applications (Theorem \ref{main tool}). The main idea in Section \ref{S: key thm} is the notion of \emph{a functional being positive relative to a set of functionals} (Definition \ref{pos}).  Then a useful characterization of this notion is given in Theorem \ref{pos char}. Some applications of this idea are given in Section \ref{S: apps}. They include a proof of the existence of numerical quadrature with positive coefficients, new proofs of Farkas type theorems,  an application to determining best approximations from certain convex cones in Hilbert space, and a specific application of the latter to determine best approximations that are also \emph{shape-preserving}.   Finally, in Section \ref{S: van char}, we note that the notion of a functional \emph{vanishing relative to a set of functionals} has a similar characterization. 

\section{The Key Theorem}\label{S: key thm} 

  The  classical Hahn-Banach separation theorem (see, e.g., \cite[p. 417]{dusc;58}) may be stated as follows. (We shall restrict our attention throughout this paper to \emph{real} linear spaces although the general results have analogous versions in complex spaces as well.) 

\begin{thm}\label{big thm}{\rm (Separation Theorem)} If $K_1$ and $K_2$ are disjoint closed convex subsets of a {\rm(real)} locally convex linear topological space $L$, and $K_1$ is compact, then there exists a continuous linear functional $f$ on $L$ such that \be\label{big thm: eq 1} \sup_{x\in K_2} f(x) < \inf_{y\in K_1} f(y). \ee \end{thm} 

The main tool of this paper (Theorem \ref{main tool}) is the special case of Theorem \ref{big thm} when $K_1$ is a single point,  $L=X^*$ is the dual space of the Banach space $X$, and $X^*$ is endowed with the  weak$^*$ topology. In the latter case, the weak*  continuous linear functionals on $X^*$ are precisely those of the form $\hat{x}$, for each $x\in X$, defined on $X^*$ by \be\label{w^*: eq 1} \hat{x}(x^*):=x^*(x)  \mbox{ \quad for each $x^* \in X^*$} \ee  (see, e.g., \cite[p. 422]{dusc;58}). It is well-known that $X$ and  $\widehat{X}:=\{ \hat{x} \mid x\in X\} \subset X^{**}$ are isometrically isomorphic. $X$ is called reflexive if $\widehat{X}=X^{**}$. 

Thus the main tool of this paper can be stated as the following corollary of Theorem \ref{big thm}. 

\begin{thm}\label{main tool}{\rm (Main Tool)} Let $X$ be a  (real) normed linear space, $\Gamma$ a weak* closed convex cone in $X^*$, and $x^* \in X^* \setminus \Gamma$.  Then there exists $x \in X$ such that  \be\label{eq 1: main tool} \sup_{y^* \in \Gamma}  y^*(x)=0 < x^*(x). \ee \end{thm}  

\pf By Theorem \ref{big thm} (with $L=X^*$ endowed with its weak* topology, $K_1=\Gamma$, and $K_2=x^*$), we deduce that there exists $x\in X$ such that  \be \label{eq 2: main tool} \sup_{y^* \in \Gamma} y^*(x) <  x^*(x). \ee Since $\Gamma$ is a cone, $0\in\Gamma$ so that $\sup_{y^* \in \Gamma} y^*(x)\ge 0$.  But if $\sup_{y^* \in \Gamma} y^*(x)> 0$, then there exists $y_0^*\in \Gamma$ such that $ y_0^*(x)> 0$. Since $\Gamma$ is a cone, $ny_0^*\in \Gamma$ for each $n\in \N$ and hence $\lim_{n \to \infty} ny_0^*(x)=\infty$. But this contradicts the fact that this expression is bounded above by $ x^*(x)$ from the inequality (\ref{eq 2: main tool}). This proves (\ref{eq 1: main tool}).  \qed 

\medskip  

Recall that the \textbf{dual cone (annihilator) in $X^*$ of a set $S \subset X$}, denoted $S^\ominus$ ($S^\perp$), is defined by  \[S^\ominus:=\{ x^* \in X^* \mid  x^*(s) \le 0 \mbox{\quad for each $s \in S$}\}\]  \[ (S^\perp:=S^\ominus\cap[-S^\ominus]=\{ x^* \in X^* \mid x^*(s)=0 \mbox{\quad for each $s \in S$} \}).  \] Clearly, $S^\ominus$ ($S^\perp$)  is a weak$^*$ closed convex cone (subspace) in $X^*$. Similarly, if $\Gamma \subset X^*$, then the \textbf{dual cone (annihilator) in $X$ of $\Gamma$}, denoted $\Gamma_\ominus$ ($\Gamma_\perp$),  is  defined by $$ \Gamma_\ominus:=\{ x \in X \mid x^*(x) \le 0  \mbox{  for all $x^* \in \Gamma$} \}$$ \[ ( \Gamma_\perp:=\Gamma_\ominus\cap[-\Gamma_\ominus]=\{ x \in X \mid x^*(x) =0  \mbox{  for all $x^* \in \Gamma$} \}). \] Clearly, $\Gamma_\ominus$  ($\Gamma_\perp$) is a closed convex cone (subspace) in $X$.   The \textbf{conical hull} of a set $S \subset X$, denoted $\cone(S)$, is the smallest convex cone that contains $S$, i.e., the intersection of all convex cones that contain $S$. Equivalently,  \be\label{cone} \cone(S):=\left\{  \sum_1^n \rho_i s_i \mid  \rho_i \ge 0, s_i \in S, n< \infty \right\}. \ee The (norm) closure of $\cone(S)$ will be denoted by $\ocone(S)$. If $S \subset X^*$, then the weak$^*$ closure of $\cone(S)$ will be denoted by $w^*\!\!-cl (\cone(S))$. 

\begin{defn}\label{pos}  Let $\Gamma$ be a subset of $X^*$. An element $x^*\in X^*$ is said to be \textbf{positive relative to $\Gamma$} if $x\in X$ and $y^*(x) \ge 0$ for all $y^* \in \Gamma$ imply that $x^*(x) \ge 0$. \end{defn}    Similarly, by replacing both ``$\ge$'' signs in Definition \ref{pos} by ``$\le$'' signs, we obtain the notion of $x^*$ being \textbf{negative relative to} $\Gamma$.  The following theorem governs this situation.    \begin{thm}\label{pos char} Let $X$ be a normed linear space, $\Gamma \subset X^*$, and $x^*\in X^*$.  Then the following statements are equivalent:  \begin{enumerate}  \item[{\rm(1)}] $x^*$ is positive relative to $\Gamma$.  \item[{\rm(2)}] $x^*$ is negative relative to $\Gamma$.  \item[{\rm(3)}] $\Gamma_\ominus \subset (x^*)_\ominus$.  \item[{\rm(4)}] $x^* \in  w^*\!\!-cl (\cone(\Gamma))$,  the weak$^*$ closed conical hull of $\Gamma$.  \end{enumerate}  Moreover, if $X$ is reflexive, then each of these statements is equivalent to   \medskip  \indent \!\!\!\!\!{\rm (5)} $x^* \in \ocone{\Gamma}$, the {\rm(}norm{\rm)} closed conical hull of $\Gamma$.  \end{thm}    \pf $(1) \Rightarrow (2)$. Suppose (1) holds, $z\in X$, and $y^*(z) \le 0$ for all $y^* \in \Gamma$. Then $y^*(-z) \ge 0$ for all $y^* \in \Gamma$. By (1), $x^*(-z) \ge 0$ or $x^*(z)\le 0$. Thus (2) holds.    $(2) \Leftrightarrow (3)$. Suppose (2) holds. If $x \in \Gamma_\ominus$, then $y^*(x) \le 0$ for all $y^* \in \Gamma$. By (2), $x^*(x) \le 0$ so $x \in (x^*)_\ominus$. Thus (3) holds. Conversely, if (3) holds, then (2) clearly holds.    $(3) \Rightarrow (4)$.  If (4) fails, then $x^* \notin w^*\!\!-cl (\cone(\Gamma))$. By Theorem \ref{main tool}, there exists $x \in X$ such that    \be\label{pos char: eq 1}  \sup\{ y^*(x) \mid y^* \in \cone(\Gamma) \}=0<x^*(x).  \ee  In particular, $y^*(x) \le 0$ for all $y^*\in \Gamma$, but $x^*(x)>0$.  Thus $x^*$ is not negative relative to $\Gamma$. That is, (2) fails.    $(4) \Rightarrow (1)$. If (4) holds, then there is a net $(y_\alp^*) \in \cone(\Gamma)$ such that $x^*(x) =\lim_\alp y_\alp^*(x)$ for all $x\in X$. If $z \in X $ and $y^*(z) \ge 0$ for all $y^* \in \Gamma$, then, in particular, $y_\alp^*(z) \ge 0$ for all $\alp$ implies that $x^*(z)=\lim_\alp y_\alp^*(z) \ge 0$. That is, $x^*$ is positive relative to $\Gamma$. Hence (1) holds, and the first four statements are equivalent.    Finally, suppose that $X$ is reflexive.  It suffices to show that   $\ocone(\Gamma)=w^*\!\!-cl(\cone(\Gamma)$. Since $X$ is reflexive,  the weak topology and the weak$^*$ topology agree on $X^*$ (see, e.g., \cite[Proposition 3.113]{fhhmz;11}). But a result of Mazur (see, e.g., \cite[Theorem 3.45]{fhhmz;11}) implies that a convex set is weakly closed if and only if it is norm closed. \qed      \medskip  In a Hilbert space $H$, we denote the inner product of $x$ and $y$ by $\la x, y\ra$ and the norm of $x$ by $\|x\|=\sqrt{\la x, x\ra}$. Then, owing to the Riesz Representation Theorem which allows one to identify $H^*$ with $H$,  Definition \ref{pos} may be restated as follows.    \begin{defn} A vector $x$ in a Hilbert space $H$ is said to be \textbf{positive relative to the set $\Gamma \subset H$} if $y\in H$ and $\la z, y \ra \ge 0$ for all $z \in \Gamma$   imply that $\la x, y\ra \ge 0$.  \end{defn}    Similarly,  in a Hilbert space $H$, we need only one notion of a dual cone (annihilator).  Namely,  if $S \subset H$, then   \[  S^\ominus: =\{ x \in H \mid \la x, y \ra \le 0 \mbox{  for all $y \in S$} \}  \]  \[ \left(S^\perp=S^\ominus \cap (-S^\ominus)=\{x\in H \mid \la x, y\ra =0  \mbox{  for all $y\in S$}\} \right).  \]  Since a Hilbert space is reflexive, we obtain the following immediate consequence of Theorem \ref{pos char}.    \begin{cor}\label{H pos char}  Let $H$ be a Hilbert space, $\Gamma \subset H$, and $x\in H$.  Then the following statements are equivalent:  \begin{enumerate}  \item[{\rm (1)}] $x$ is positive relative to $\Gamma$. \item[{\rm (2)}]  $x$ is negative relative to $\Gamma$.     \item[{\rm (3)}] $\Gamma^\ominus \subset (x)^\ominus$. \item[{\rm (4)}] $x \in \ocone{\Gamma}$. \end{enumerate} \end{cor} 

\medskip Well-known examples of reflexive spaces are  finite-dimensional spaces, Hilbert spaces, and the $L_p$ spaces for $1< p < \infty$.  (The spaces $L_1$ and $C(T)$, for $T$ compact, are never reflexive unless they are finite-dimensional.) 

For a general convex set, we have the following relationship. 

\begin{lem}\label{convex} Let $K$ be a convex subset of a normed linear space $X$. Then \be\label{convex: eq 1} (K^\ominus)_\ominus=\ocone(K) \ee \end{lem} 

\pf  By definition, $K^\ominus=\{ x^*\in X^* \mid x^*(K) \le 0\}$. Thus \begin{eqnarray*} (K^\ominus)_\ominus&=&\{ x \in X \mid x^*(x) \le 0  \mbox{  for all $x^* \in K^\ominus$} \}\\   &=&\{x\in X \mid x^*(x) \le 0   \mbox{   for each $x^*$ such that $x^*(K)\le 0$} \}\\   &\supset& K.   \end{eqnarray*}   Since $(K^\ominus)_\ominus$ is a closed convex cone, it follows that             $(K^\ominus)_\ominus \supset \ocone(K)$. If the lemma were false, then there would exist $x \in (K^\ominus)_\ominus \setminus \ocone(K)$.  By Theorem \ref{big thm}, there exists $x^* \in X^*$ such that $\sup x^*[\ocone(K)] < x^*(x)$.  Arguing as in the proof of Theorem \ref{main tool}, we deduce that    $ \sup x^*[\ocone(K)]=0 < x^*(x)$. But this contradicts the fact that $x^* \in K^\ominus$ and $x\in (K^{\ominus})_\ominus$.  \qed  

 \medskip  \begin{cor}\label{cone cor} If $C$ is a nonempty  subset of $X$, then $C$ is a closed convex cone in $X$ if and only if   \be\label{cone cor: eq 1}  C=(C^\ominus)_\ominus .  \ee  \end{cor}      

\medskip It follows  that \emph{every} closed convex cone has the same special form. More precisely, we have the following easy consequence. 

\begin{lem}\label{char cone} Let $X$ be a normed linear space and let $C$ be a nonempty subset of $X $. Then the following statements are equivalent: \begin{enumerate} \item[{\rm (1)}] $C$ is a closed convex cone. \item[{\rm (2)}] There exists a set $\Gamma \subset X^*$ such that \[ C:=\{ x\in X \mid y^*(x) \le 0  \mbox{  for each $y^*\in \Gamma$} \}. \] {\rm(}In fact, $\Gamma=C^\ominus$ works.{\rm)} \item[{\rm (3)}] There exists a set $\widetilde{\Gamma} \subset X^*$ such that \[ C:=\{ x\in X \mid y^*(x) \ge 0  \mbox{  for each $y^*\in \widetilde{\Gamma}$} \}. \] {\rm(}In fact, $\widetilde{\Gamma}=-C^\ominus$ works.{\rm)}  \end{enumerate}  \end{lem} 

We will need the following fact that goes back to Minkowski (see, e.g., \cite[Lemma 6.33]{deu;01}).   

\begin{fact}\label{mink} If a nonzero vector $x$ is a positive linear combination of the vectors $x_1, x_2, \dots, x_n$, then $x$ is a positive linear combination of a linearly independent subset of $\{x_1, x_2, \dots, x_n\}$. \end{fact}    \begin{thm}\label{cone} Let $X$ be a reflexive Banach space and $\Gamma \subset X^*$ be weakly compact. Suppose there exists $y \in X$ such that $y^*(y)> 0$ for each $y^*\in \Gamma$ and $\dim(\Gamma)=n$ {\rm (}so  $\Gamma$ contains a maximal set of $n$ linearly independent vectors{\rm)}.  Then each nonzero $x^* \in \ocone(\Gamma)$ has a representation as  \be\label{cone: eq 1}  x^*=\sum_1^m \rho_i y_i^* ,  \ee  where $m \le n$, $\rho_i > 0$ for $i=1, 2, \dots, m$, and $\{ y_1^*, y_2^*, \dots, y_m^*\}$ is a linearly independent subset of $\Gamma$.  \end{thm}    \pf Let $\del:=\inf\{ y^*(y) \mid y^* \in \Gamma \}$. If $\del=0$, then there exists a sequence $(y_n^*)$ in $\Gamma$ such that $\lim y_n^*(y)=0$.  By the Eberlein-\u{S}mulian Theorem (see, e.g., \cite[p. 129]{fhhmz;11}), $\Gamma$ is weakly sequentially compact, so  there is a subsequence $(y_{n_k}^*)$  which converges weakly to $y^* \in \Gamma$ and, in particular, $0=\lim y^*_{n_k}(y)=y^*(y)>0$,  which is absurd. Thus $\del >0$.    Let $x^* \in \ocone(\Gamma)\setminus \{0\}$. Then there exists a sequence $(x_N^*)_1^\infty$ in $\cone(\Gamma)$ such that $\|x_N^* - x^*\| \to 0$. Since $x^* \neq 0$, we may assume that $x_N^*\neq 0$ for all $N$. Then $(x_N^*)$ is bounded, say $c:=\sup_N \|x_N^*\|  < \infty$, and   \be\label{cone: eq 2}  x_N^*=\sum_{i \in F_N} \rho_{N,i} x_{N, i}^*  \ee   for some scalars $\rho_{N,i}\ge 0$, $x_{N,i}^* $ in $\Gamma$, and $F_N$ is finite. By the hypothesis $\dim(\Gamma)=n$ and Fact \ref{mink}, we may assume that $F_N =\{1,2, \dots, n\}$. Thus we have that   \be\label{cone: eq 3}   x_N^*=\sum_1^n \rho_{N,i} x_{N,i}^*   \ee   where $\rho_{N,i} \ge 0$ for all $i$. Now   \be\label{cone: eq 4}   x_N^*(y)=\sum_{i=1}^N \rho_{N, i} x^*_{N, i}(y) \ge \sum_{i=1}^n \rho_{N, i} \del \ge \rho_{N, i} \del  \mbox{  for each $i=1,2, \dots, n$. }   \ee  Thus, for each $i \in \{1, \dots, n\}$, we have  \be\label{cone: eq 5}  \rho_{N, i}\le(\del)^{-1}x_N^*(y) \le (\del)^{-1}\|x_N^*\| \|y\| \le(\del)^{-1}\|y\|c < \infty.  \ee   This shows that, for each $i=1, 2, \dots, n$, the sequence of scalars $(\rho_{N,i})$ is bounded. By passing to a subsequence, we may assume that there exist $\rho_i \ge 0$ such that $\rho_{N, i} \to \rho_i$ for each $i$.      Since $\Gamma$ is weakly sequentially compact, by passing to a further subsequence, say $(N')$ of $(N)$, we may assume that for each $i=1, 2, \dots, n$, there exist $y_i^* \in    \Gamma$ such that $x_{N', i}^* \to y_i^*$ weakly. Thus, for all $x\in X$, we have   \be\label{cone: eq 6}   x^*(x)=\lim_{N'} \sum_{i=1}^n \rho_{N', i} x_{N', i}^*(x)=\sum_{i=1}^n \rho_i y_i^*(x).   \ee   That is, $x^*=\sum_1^n \rho_i y_i^*$. By appealing to Fact \ref{mink}, we get the representation (\ref{cone: eq 1}) for $x^*$.   \qed    \medskip  Again, in the case of Hilbert space, this result reduces to the following fact that was first established by Tchakaloff \cite{tch;57}, who used it to prove the existence of quadrature rules having positive coefficients (see also Theorem \ref{quad} below).     \begin{cor}\label{coneh} Let $H$ be a Hilbert space and $\Gamma \subset H$ be weakly compact. Suppose there exists $e \in H$ such that $\la y, e \ra > 0$ for each $y\in \Gamma$ and $\dim(\Gamma)=n$ {\rm (}so  $\Gamma$ contains a maximal set of $n$ linearly independent vectors{\rm)}.  Then each nonzero $x \in \ocone(\Gamma)$ has a representation as  \be\label{coneh: eq 1}  x=\sum_1^m \rho_i y_i ,  \ee  where $m \le n$, $\rho_i > 0$ for $i=1, 2, \dots, m$, and $\{ y_1, y_2, \dots, y_m\}$ is a linearly independent subset of $\Gamma$.  \end{cor}     

\section{Some Applications of Theorem \ref{pos char}}\label{S: apps} 

In this section we show the usefulness of Theorem \ref{pos char} by exhibiting a variety of different applications.  

\subsection{An Application to the Existence of Positive \newline Quadrature Rules} 

In the first application, we show the existence of quadrature rules that are exact for polynomials of degree at most $n$, are based on a set of  $n+1$ points, and have positive coefficients. Let $\mathcal{P}_n$ denote the set of polynomials of degree (at most) $n$ regarded as a subspace of $C[a, b]$. That is, $\cP_n$ is endowed with the norm $\|x\|=\max\{ |x(t)| \mid a \le t \le b\}$. Define the linear functionals $x^*$ and $x_t^*$ on $X:=\cP_n$ by \be\label{integral1} x^*(x):=\int_a^b x(t)\,dt   \mbox{\quad     for all  $x\in X$   } \ee and  \be\label{integral2} x^*_t(x):=x(t)     \mbox{\quad     for all  $x\in X$.   } \ee  

\begin{thm}\label{quad}{\rm (Numerical Quadrature)} Let $X=\mathcal{P}_n$. Then there exists $m \le n+1$ points $a\le t_1 < t_2 < \cdots < t_m \le b$ and $m$  scalars $w_i > 0$ such that   $x^*=\sum_1^m w_i x^*_{t_i}$. More explicitly, \be\label{quad: eq 1} \int_a^b x(t)\,dt =\sum_1^m w_i x(t_i)   \mbox{   for all $x\in X$.  } \ee \end{thm} 

\pf  First note that  $x^*$ is positive relative to the set $\Gamma:=\{ x_t^* \mid t\in [a,b] \}$, since a function that is nonnegative at each point in $[a,b]$ must have a nonnegative integral.  Since $X$ is finite-dimensional, it is reflexive. By Theorem \ref{pos char}(4), we have $x^*\in \ocone(\Gamma)$. Next note that for the identically 1 function $e$ on $[a,b]$, we have $x_t^*(e)=1$ for all $t\in [a,b]$.   Further, it is easy to check that $\Gamma$ is a closed and bounded subset of $X^*$, hence is compact since in a finite-dimensional  space $X$ all linear vector space topologies on $X^*$ coincide (see, e.g., \cite[Corollary 3.15]{fhhmz;11}). Finally, since $\dim X^*=\dim X=n+1$, we can apply Theorem \ref{cone}  to get the   result.  \qed 

\subsection{Applications Related to Farkas Type Results} 

In this section we note that the so-called  Farkas Lemma is a consequence of Theorem \ref{pos char}.  According to Wikipedia,  \begin{quote} Farkas' lemma is a solvability theorem for a finite system of linear inequalities in mathematics. It was originally proven by the Hungarian mathematician Gyula Farkas \cite{far;02}. Farkas' lemma is the key result underpinning the linear programming duality and has played a central role in the development of mathematical optimization (alternatively, mathematical programming ). It is used amongst other things in the proof of the  Karush-Kuhn-Tucker theorem in nonlinear programming. \end{quote} 

  Since the setting for this result is in a Hilbert space, we will be appealing to the Hilbert space version of Theorem \ref{pos char}, namely, Theorem \ref{H pos char}. 

\begin{thm}\label{F1} Let $H$ be a Hilbert space and $\{b, a_1, a_2, \dots, a_m\} \subset H$. Then exactly one of the following two systems has a solution: \par System 1: $\sum_1^m y_ia_i=b$ for some $y_i \ge 0$. \par System 2: There exists $x\in H$ such that $\la a_i, x\ra \le 0$ for $i=1, \dots, m$ and $\la b, x \ra >0$. \end{thm} 

\pf Letting $\Gamma:=\{a_1, a_2, \dots, a_m\}$, we see that the cone generated by $\Gamma$ is finitely generated and, as is well-known,  must be closed (see, e.g., \cite[Theorem 6.34]{deu;01}). Hence $\ocone(\Gamma)=\{\sum_1^m\rho_ia_i \mid \rho_i \ge 0\}$.  

Clearly,  system 1 has a solution  if and only if  $b\in \ocone\{a_1, a_2, \dots, a_m\}$. By Theorem \ref{H pos char}, system 1 has a solution if and only if $b$ is negative relative to $\Gamma:=\{a_1, a_2, \dots, a_m\}$.   

But obviously,  system 2 has a solution if and only if $b$ is not negative relative to $\Gamma$.  This completes the proof.  \qed 

If this theorem is given in its (obviously equivalent) matrix formulation, then it can be stated as in the following theorem. This is the version given by Gale, Kuhn, and Tucker \cite{gkt;51} (where vector inequalities are interpreted componentwise). 

\begin{thm}\label{F2} Let $A$ be an $m \times n$ matrix and $b\in \R^n$.  Then exactly one of the following two systems has a solution: \par System 1: $A^Ty=b$ and $y\ge 0$ for some $y\in \R^m$. \par System 2:  There exists $x\in R^n$ such that $Ax \le 0$ and $\la b, x\ra >0$. \end{thm} 

The next theorem extends a result of Hiriart-Urruty and Lemar\'echal \cite[Theorem 4.3.4]{hul;93} who called it a \emph{generalized Farkas theorem}. 

\begin{thm}\label{F3} Let $J$ be a index set and $(b, r)$ and $(s_j, p_j) \in \R^n \times \R$ for all $j\in J$.  Suppose that the system of inequalities  \be\label{eq 1: F3}  \la s_j, x \ra \le p_j  \ee \par  has a solution $x\in \R^n$. Then the following statements are equivalent: \begin{enumerate} \item[{\rm(1)}] $\la b, x\ra \le r$ for all $x$ that satisfy relation $(\ref{eq 1: F3})$. \item[{\rm (2)}] $( b, r) \in \ocone\{(s_j, p_j) \mid j \in J \}$. \item[{\rm (3)}]  $(b,r) \in \ocone \left(\{(0,1)\} \cup \{ (s_j, p_j)\mid j \in J \}\right)$. \end{enumerate} \end{thm} 

\pf First note that $x \in \R^n$ is a solution to  (\ref{eq 1: F3}) if and only if  $(x, -1) \in \R^n \times \R$ is a solution to  \be\label{eq 2: F3} \la (s_j, p_j), (x, -1) \ra \le 0  \mbox{   for all $j\in J$}. \ee Using this fact, we see that statement (1) holds $\iff$  \be\label{eq 3: F3} \la b,  x \ra - r \le 0  \mbox{   for all $x$ that satisfy (\ref{eq 2: F3})} \ee $\iff$ \be\notag \la (b, r), (x,-1) \ra \le 0   \mbox{  for all $x$ such that } \la (s_j, p_j), (x, -1) \ra \le 0 \mbox{  for all $j\in J$.} \ee But the last statement just means that $(b,r)$ is negative relative to the set $\{(s_j, p_j) \mid j \in J\}$.  By Theorem \ref{H pos char}, this is equivalent to statement (2).  Thus we have proved $(1) \iff (2)$.  

Clearly (2) implies (3) since the conical hull in (3) is larger than that of (2).  Finally, the same proof of (3) implies (1) as given in \cite[Theorem 4.3.4]{hul;93} works here.  \qed 

\begin{rem} Hirriart-Urruty and Lemar\'echal \cite[Theorem 4.3.4]{hul;93} actually proved the equivalence of statements (1) and (3) of Theorem \ref{F3}.  The sharper equivalence of statements (1) and (2) proven above was seen to be a simple consequence of Corollary \ref{H pos char}. \end{rem} 

\subsection{An Application to Best Approximation}\label{S: Hilbert space} 

In this section we give an application to a problem of best approximation from a convex cone in a Hilbert space. We will need a special case of the following well-known  characterization of best approximations from convex sets. This characterization goes back at least to Aronszajn \cite{aro;50} in 1950 (see also \cite[Theorem 4.1]{deu;01}). The fact that every closed convex subset $C$ of a Hilbert space $H$ admits unique nearest points (best approximations) to each $x\in H$ is due to Riesz \cite{rie;34}. If $x\in H$, we denote its unique best approximation in $C$ by $P_C(x)$. 

\begin{fact}\label{aron} Let $C$ be a closed convex set in a Hilbert space $H$, $x\in H$, and $x_0\in C$.  Then $x_0=P_C(x)$ if and only if $x-x_0\in (C-x_0)^\ominus$, i.e., \be\label{aron; eq 1} \la x-x_0, y-x_0 \ra \le 0  \mbox{  for each $y\in C$}. \ee \end{fact} 

In the special case when $C$ is a closed convex cone, Moreau \cite{jjm;62} showed, among other things,  that this result could be sharpened to the following. 

\begin{fact}\label{moreau} Let $C$ be a closed convex cone in the Hilbert space $H$, $x\in H$, and $x_0 \in C$.  Then $x_0=P_C(x)$ if and only if $x-x_0\in C^\ominus \cap x_0^\perp$, i.e.,  \be\label{eq 1: mor} \la x-x_0, y\ra \le 0  \mbox{  for all $y\in C$ and }   \la x-x_0, x_0 \ra =0. \ee Moreover, $H=C \boxplus C^\ominus$, which means that each $x\in H$ has a unique representation as $x=c+c'$ where $c\in C$, $c'\in C^\ominus$, and $\la c, c' \ra=0$. In fact,  \be\label{eq 2: mor} x=P_C(x) + P_{C^\ominus}(x)  \mbox{  for each $x \in H$}. \ee \end{fact} 

(For  proofs of these facts, see, e.g., \cite[Theorems 4.1, 4.7, and 5.9]{deu;01}.)

\setlength{\abovecaptionskip}{12pt} \setlength{\belowcaptionskip}{0pt} \setlength{\textfloatsep}{-5pt} \setlength{\intextsep}{-5pt} \begin{wrapfigure}{r}{0.45\textwidth}   \vspace{12pt}   \centering   \includegraphics*[height=1.75in]{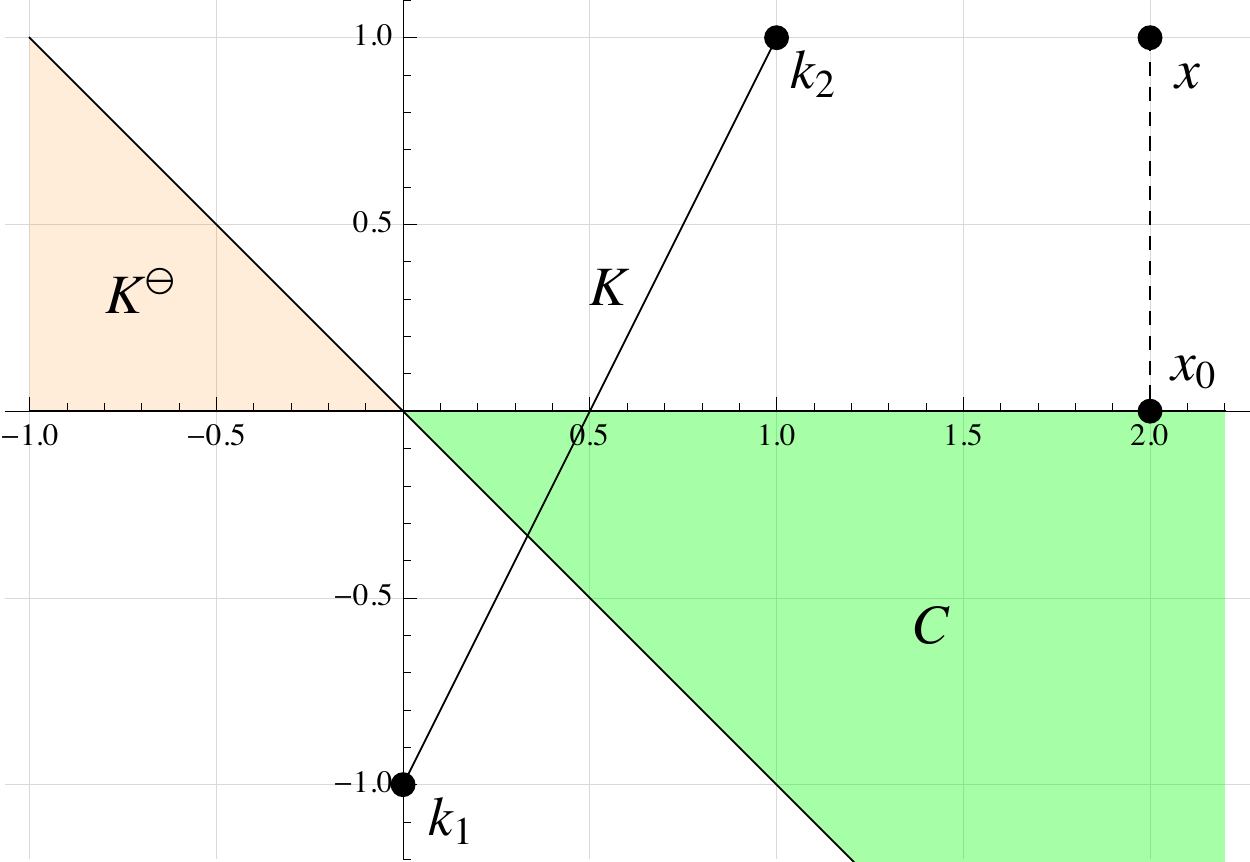}   \caption{\label{fig:ExampleCone}} \end{wrapfigure} The main result of this section is Theorem \ref{deu thm}.  To provide some motivation, we  exhibit a simple example. 

\begin{ex}\label{simple ex} Let $X=\ell_2(2)$ denote Euclidean 2-space, $K$ denote the line segment joining the two points $k_1=(0,-1)$ and $k_2=(1,1)$, and $C=-K^\ominus$,  i.e., $C=\{y\in X \mid  \la y, k_i \ra \ge 0  \mbox{   for $i=1, 2$} \}$, see Figure~\ref{fig:ExampleCone}. Let $x=(2,1)$ and $x_0=(2, 0)$. Then $x_0=x+k_1=P_C(x)$. \end{ex} 

\begin{thm}\label{deu thm} Let $K$ be a compact set in the  Hilbert space $H$, and suppose that there exists $e \in H$  such that  \begin{enumerate} \item[{\rm(i)}] $\la k, e\ra >0$ for all $k\in K$, and \item[{\rm(ii)}] $\dim K=n$.  \end{enumerate} Let $C$ be {\rm(}the closed convex cone{\rm)} defined by \be\label{deu thm: eq 1} C:=-K^\ominus=\{ y \in H \mid \la y, k\ra \ge 0  \mbox{    for all $k \in K$} \}. \ee Let $x \in H \setminus C$ and $x_0\in C$. Then the following statements are equivalent: \begin{enumerate} \item[{\rm(1)}] $x_0=P_C(x)$; \item[{\rm(2)}]  \be\label{deu thm: eq 2} x_0=x+ \sum_1^m \rho_i k_i, \ee where $1\le m \le n$, $\rho_i >0$,  $k_i \in K$ for $i=1, 2, \dots, m$, where the vectors $ k_1, k_2, \dots, k_m$ are linearly independent, and $\la k_i, x_0 \ra =0$ for $i=1, 2, \dots, m$. \end{enumerate} 

Moreover, if $\dim H=n$ and $x_0\ne 0$ in any of the two statements,  then $m\le n-1$. \end{thm} 

\pf $(1) \Rightarrow (2)$. If (1) holds, then by Fact \ref{aron}, we have that $\la x-x_0, y\ra \le 0$ for all $y \in C$ and $\la x-x_0, x_0\ra =0$.  Thus if $y \in X$ and $\la y, k\ra \ge 0$ for all $k\in K$, then $y\in C$ so that $-\la x_0-x, y \ra =\la x -x_0, y \ra \le 0$, so $\la x_0 -x, y \ra \ge 0$.  Thus $x_0-x$ is positive relative to $K$. By Corollary \ref{H pos char}, we see that $x_0-x \in \ocone(K)$. By Lemma  \ref{coneh}, we have that $x_0-x=\sum_1^m \rho_i k_i$, where $\rho_i >0 $ for all $i$, $m \le n$, and the set $\{ k_1, \dots, k_m\}$ is linearly independent. Also, since $\la x-x_0, x_0 \ra =0$, we see that  \be\label{deu thm: eq 1} \sum_1^m \rho_i\la k_i, x_0\ra= \left\la \sum_1^m \rho_i k_i, x_0 \right\ra =\la x_0-x, x_0\ra=0 \ee which, since $\la k_i,x_0 \ra \ge 0$ and $\rho_i > 0$ for all $i$, implies that $\la k_i, x_0 \ra =0$ for all $i$.  Thus (2) holds. 

$(2) \Rightarrow (1)$.  If (2) holds, then $x_0-x=\sum_1^m \rho_i k_i$ and $\la k_i, x_0\ra =0$ for all $i$.  Thus for all $y\in C$ we have $\la x_0-x, y\ra =\sum_1^m \la k_i, y\ra \ge 0$ and $\la x_0-x, x_0\ra =\sum_1^m \rho_i \la k_i, x_0 \ra =0$. In other words, $x-x_0 \in C^\circ \cap x_0^\perp$.  By Fact \ref{moreau}, we see that $x_0=P_C(x)$, i.e., (1) holds. This proves the equivalence of the two statements. 

Finally, if $\dim X =n$ and $x_0 \ne 0$ in any of the two statements,  then we see that  $\la x_0, k_i\ra =0$  for each $i=1, \dots, m$. But the null space of $x_0$, i.e., $x_0^\perp:=\{y \in X \mid \la x_0, y\ra =0 \}$, is an $(n-1)$-dimensional subspace of the  $n$-dimensional space $X$. Since $\{k_1, \dots, k_m\}$ is a linearly independent set contained in $x_0^\perp$, we must have $m\le n-1$.  \qed  

\begin{rem} It is worth noting that  if either of the equivalent statements (1) or (2) holds in Theorem \ref{deu thm}, then there exists at least one $i$ such that $\la x, k_i\ra < 0$. \end{rem}  

To see this, assume (2) holds. Then \begin{eqnarray*} \sum_1^m \rho_i\la k_i, x\ra &=& \left\la \sum_1^m \rho_ik_k, x\right\ra =\la x_0-x, x \ra=\la x_0-x, x-x_0 \ra\\ &=&-\|x_0-x\|^2 < 0  \end{eqnarray*} which, since $\rho_i > 0$ for each $i$, implies that $\la k_i, x\ra < 0$ for some $i$.    

The following corollary of Theorem \ref{deu thm} shows that in certain cases, one can even obtain an \emph{explicit formula} for the best approximation to any vector. 

\begin{cor}\label{cor of deu thm}Let $K=\{k_1, k_2, \dots, k_n\}$ be an orthonormal subset of the Hilbert space $H$, and suppose that there exists $e\in H$ such that $\la e, k_i\ra > 0$ for each $i$. Let  \be C:=-K^\ominus=\{y \in H \mid \la y, k_i \ra \ge 0 \mbox{  for each $i$}\}, \ee and $x\in H$.  Then \be\label{eq 1: cor of deu thm} P_C(x)=x+\sum_{i=1}^n \max\{ 0, -\la x, k_i \ra\} k_i. \ee \end{cor} 

\pf If $x\in C$, then $P_C(x)=x$ and $\la x, k_i \ra \ge 0$ for each $i$ implies that $\max\{0, -\la x, k_i \ra \}=0$ for each $i$ and thus formula (\ref{eq 1: cor of deu thm}) is correct. Hence we may assume that $x\in H \setminus C$. 

Let $x_0=x + \sum_{i=1}^n \rho_i k_i$, where $\rho_i=\max\{0, -\la x, k_i\ra \}$. It suffices to show that $x_0=P_C(x)$.  Let $J=\{ j \in \{1,2, \dots, n\} \mid \la x, k_j \ra < 0\}$.  Since $x\notin C$, we see that $J$ is not empty, $\rho_j=-\la x, k_j\ra$ for each $j\in J$,  $\rho_i =0$ for all $i \notin J$, and $x_0=x+\sum_{j\in J} \rho_j k_j$.  

Using the orthonormality of the set $K$, we see that for all $j \in J$,  \be\label{eq 1: cor deu thm} \la x_0, k_j \ra= \la x, k_j \ra +\left\la\sum_{i \in J} \rho_i k_i, k_j\right\ra=\la x, k_j\ra +\rho_j=0, \ee and for each $i \notin J$, \be\label{eq 2: cor deu thm} \la x_0, k_i \ra=\la x, k_i\ra + \sum_{j \in J}\rho_i \la k_j, k_i \ra=\la x, k_i \ra \ge 0. \ee The relations (\ref{eq 1: cor deu thm}) and (\ref{eq 2: cor deu thm}) together show that $x_0 \in C$. Finally, the equality (\ref{eq 1: cor deu thm}) shows that Theorem \ref{deu thm}(2) is verified.  Thus $x_0=P_C(x)$, and  the proof is complete.  \qed 

We next consider \emph{two} alternate versions of Theorem \ref{deu thm} which may be more useful for the actual computation of best approximations from finitely generated convex cones. 

We  consider the following scenario. Let $H$ be a Hilbert space, $\{k_1, k_2, \dots, k_m\}$ a finite subset of $H$, and $C$ the convex cone generated by $K$:  \be\label{eq 1: gdeu thm} C:=\left\{\sum_{i=1}^m \rho_ik_i \biggm| \mbox{  $\rho_i \ge 0$ for all $i$} \right\}=\ocone\{k_1, k_2, \dots, k_m\}. \ee 

By definition of the dual cone, we have \begin{eqnarray}\label{eq 2: gdeu thm} C^\ominus&=&\{ y\in H \mid \la y, c \ra \le 0  \mbox{ for all $c\in C$}\}\notag \\   &=&\{ y \in H \mid \la y, k_i  \ra \le 0 \mbox{  for all $i=1, \dots, m $} \}. \end{eqnarray}

As an easy consequence of a theorem of the first author characterizing best approximations from a polyhedron  (\cite[Theorem 6.41]{deu;01}), we obtain the following. 

\begin{thm}\label{deu poly} Let $\{k_1, k_2, \dots, k_m\}$ be a finite subset of the Hilbert space $H$, and let $C$ be the finitely-generated cone defined by equation {\rm (\ref{eq 1: gdeu thm}).} Then for each $x\in H$,  \be\label{eq 1: deu poly}  P_{C^\ominus}(x) =x - \sum_1^m \rho_ik_i  \mbox{\quad  and \quad  }    P_C(x)=\sum_1^m \rho_ik_i    \ee for any set of scalars $\rho_i$ that satisfy the following three conditions: \be\label{eq 2: deu poly} \rho_i \ge 0 \qquad (i=1, 2, \dots, m) \ee \be\label{eq 3: deu poly} \la x, k_i\ra - \sum_{j=1}^m \rho_j \la k_j, k_i\ra \le 0  \qquad (i=1, 2,\dots, m) \ee and \be\label{eq 4: deu poly} \rho_i[\la x, k_i \ra- \sum_{j=1}^m \rho_j \la k_j, k_i\ra]=0 \qquad (i=1, 2, \dots, m). \ee 

Moreover, if $x\in H$ and $x_0\in C^\ominus$, then $x_0=P_{C^\ominus}(x)$ if and only if  \be \label{eq 5: deu poly} x_0=x-\sum_{i \in I(x_0)} \rho_i k_i  \mbox{  for some $\rho_i \ge 0$,  } \ee where $I(x_0): =\{ i \mid \la x_0, k_i \ra =0\}$. \end{thm} 

\pf  In \cite[Theorem 6.41]{deu;01}, take $X=H$,  $c_i=0$ and $h_i=k_i$ for all $i=1, 2, \dots, m$, and note that $Q=\{y \in H \mid  \la y, k_i \ra \le 0\}=C^\ominus$.  The conclusion of \cite[Theorem 6.41]{deu;01}  now shows that $P_{C^\ominus}(x)=P_Q(x)=x-\sum_1^m\rho_ik_i$, where the $\rho_i$ satisfy the relations (\ref{eq 2: deu poly}), (\ref{eq 3: deu poly}), and (\ref{eq 4: deu poly}). Finally, by Fact \ref{moreau}, we obtain that $P_C(x)=x-P_{C^\ominus}(x)=\sum_1^m \rho_ik_i$.   

The last statement of the theorem follows from the last statement of \cite[Theorem 6.41]{deu;01}.  \qed  

We will prove an alternate characterization  of best approximations from finitely generated cones that yields detailed information of a different kind. But first we need to recall some relevant concepts.  

For the remainder of this section, we assume that $T: H_1\to H_2$ is a bounded linear operator between the Hilbert spaces $H_1$ and $H_2$ that has closed range. Then the adjoint mapping $T^*: H_2 \to H_1$ also has closed range (see, e.g., \cite[Lemma 8.39]{deu;01}). We denote the \emph{range} and \emph{null} space of $T$ by \[ \label{eq 3: zig} \mathcal{R}(T):=\{T(x) \mid  x \in H_1\}, \qquad \mathcal{N}(T):=\{ x \in H_1 \mid  T(x)=0 \}. \] 

The following relationships between these concepts are well-known (see, e.g., \cite[Lemma 8.33]{deu;01}): \begin{eqnarray} \mathcal{N}(T)&=&\mathcal{R}(T^*)^\perp, \quad \mathcal{N}(T^*)=  \mathcal{R}(T)^\perp, \mbox{ \quad and \quad  } \\  \mathcal{N}(T)^\perp=\overline{\mathcal{R}(T^*)}&=&\mathcal{R}(T^*), \quad   \mathcal{N}(T^*)^\perp= \overline{\mathcal{R}(T^)}=\mathcal{R}(T). \label{eq 2} \end{eqnarray} 

\begin{defn} For any $y\in H_2$, the set of \textbf{generalized solutions} to the equation $T(x)=y$ is the set \[ G(y):=\{x_0 \in H_1 \mid \| T(x_0) - y\| \le \|T(x) -y\| \mbox{  for all $x\in H_1$} \}. \] \end{defn} 

Since $\mathcal{R}(T)$ is closed, it is a Chebyshev set so $G(y)$ is not empty. For each $y\in H_2$, let $T^\dagger(y)$ denote the \emph{minimal norm} element of $G(y)$.  The mapping $T^\dagger:  H_2 \to H_1$ thus defined is called the \textbf{generalized inverse} of $T$.  

The following facts are well-known (see, e.g., \cite{gro;77} or \cite[pp 177--185]{deu;01}). 

\begin{fact}\label{fact 1} \begin{enumerate} \item[{\rm (1)}] $T^\dagger$ is a bounded linear mapping. \item[{\rm (2)}] $(T^*)^\dagger=(T^\dagger)^*$. \item[{\rm (3)}] $TT^\dagger=P_{\mathcal{R}(T)}=P_{\mathcal{N}(T^*)^\perp}$. \item[{\rm (4)}] $T^\dagger T=P_{\mathcal {N}(T)^\perp}$. \item[{\rm (5)}] $TT^\dagger T=T$. \end{enumerate} \end{fact} 

As in the above Theorem \ref{deu poly}, we again let $\{ k_1, k_2, \dots, k_m\}$ be a finite subset of the Hilbert space $H$ and  $C$ be the convex cone generated by the $k_i$: \be\label{eq 1: fg cone} C=\ocone\{ k_1, k_2, \dots, k_m\}=\left\{\sum_1^m \rho_i k_i \bigm|  \rho_i \ge 0 \mbox{  for all $i$ } \right\}. \ee It follows that \begin{eqnarray} C^\ominus&=&\{y \in H \mid  \la y, c\ra \le 0  \mbox{  for all $c\in C$}\}\\ &=&\{ y \in H \mid \la y, k_i \ra \le 0  \mbox{  for all $i$}\}. \end{eqnarray} Let $S: \R^m \to H$ be the bounded linear operator defined by  \be \label{eq 2: fg cone} S(\alp) = \sum_1^m \alp_i k_i \mbox{ \qquad for all}      \quad  \alp=(\alp_1, \alp_2, \dots, \alp_m) \in \R^m. \ee If $S^*: H \to \R^m$ denotes the adjoint of $S$, then  \be \label{eq 2: zig} \la S^*(y), e_j\ra =\la y, S(e_j) \ra = \la y, k_j\ra   \mbox{   for all $j$}, \ee where $e_j$ denote the canonical bases vectors in $\R^m$, i.e., $e_j=(\del_{1j}, \del_{2j}, \dots, \del_{mj})$, and $\del_{ij}$ is Kronecker's delta---the scalar which is $1$ when $i=j$ and $0$ otherwise. 

As was noted in Fact \ref{moreau}, if $C$ is a closed convex cone in a Hilbert space $H$, then $H=C\boxplus C^\ominus$, which means that each $x \in H$ has a unique orthogonal decomposition as  $x=P_C(x)+P_{C^\ominus}(x)$.  In the case of a finitely generated cone,  we will strengthen and extend this even further by showing that   $C^\ominus$ has  an even \emph{stronger}  orthogonal decomposition as the sum of ${N}(S^*)$ and a certain subset of  $\mathcal{N}(S^*)^\perp$. For a vector $\rho=(\rho_1, \dots, \rho_m) \in \R^m$, we write $\rho \ge 0$ to mean $\rho_i \ge 0$ for each $i$. 

\begin{lem}\label{zig} The following orthogonal decomposition holds: \be\label{eq 1: zig} C^\ominus=\mathcal{N}(S^*) \boxplus \mathcal{B}, \mbox{ \quad where } \ee  \be\label{eq 2: zig} \mathcal{B}:=\{z \in H \mid z=-(S^*)^\dagger(\rho), \quad \rho\in \mathcal{N}(S)^\perp, \quad \rho \ge 0\}  \subset \; \cN(S^*)^\perp.  \ee In particular, each $c'\in C^\ominus$ has a unique representation as $c'=y+z$, where $y\in \mathcal{N}(S^*)$, $z\in \mathcal{B}$  and $\la y, z\ra=0$. \end{lem} 

\pf  We first show that $\cB \subset \cN(S^*)^\perp$.  If $z\in \cB$, then $z=-(S^*)^\dagger(\rho)$, where $\rho \in \cN(S)^\perp$. Since $\cN(S)^\perp=\cR(S^*)$ by the relation (\ref{eq 2}), we can write $z=-(S^*)^\dagger S^*(u)$, for some $u\in H$. Further, by Fact \ref{fact 1}(4), we see that $z=-P_{\cN(S^*)^\perp}(u) \in \cN(S^*)^\perp$, which proves $\cB \subset \cN(S^*)^\perp$ and thus verifies (\ref{eq 2: zig}). 

Let  \[ \cD:=\mathcal{N}(S^*) \boxplus \{z \in H \mid z=-(S^*)^\dagger(\rho), \quad \rho\in \mathcal{N}(S)^\perp, \quad \rho \ge 0\}. \]  If we can show that $\cD=C^\ominus$, then the last statement of the lemma will follow from this and relation (\ref{eq 2: zig}). Thus to complete the proof, we need to show that $\cD=C^\ominus$.   

Let $y\in C^\ominus$. For each $j=1, \dots, m$, let $\rho_j:=-\la y, k_j \ra$. Since $y\in C^\ominus$, it follows that  $\rho_j \ge 0$ and so $\rho \ge 0$. To see that $\rho \in \mathcal{N}(S)^\perp$, take any $\eta \in \mathcal{N}(S)$.  Then $S(\eta)=0$ and \begin{eqnarray} \la \rho, \eta \ra&=&\sum_1^m \rho_j \eta_j=-\sum_1^m\la y, k_j\ra \eta_j \\ &=&-\la y, \sum_1^m \eta_j k_j \ra=-\la y, S(\eta) \ra=0. \end{eqnarray} Since $\eta \in \mathcal{N}(S)$ was arbitrary, it follows that $\rho \in \mathcal{N}(S)^\perp$. 

 The definition of $\rho_j$ yields \[ \la \rho, e_j \ra=\rho_j=-\la y, k_j\ra =-\la y,S(e_j) \ra =-\la S^*(y), e_j \ra=\la -S^*(y), e_j\ra. \] Since this holds for all the basis vectors $e_j$, it follows that $\rho=-S^*(y)$. Further, by Fact \ref{fact 1}(4), we see that $(S^*)^\dagger S^*=P_{\mathcal{N}(S^*)^\perp}$ and hence we can write \[ y=y-(S^*)^\dagger S^*(y) +(S^*)^\dagger S^*(y)=y_0 -(S^*)^\dagger(\rho), \]  where  \[ y_0:=[I-(S^*)^\dagger S^*](y)=[I-P_{\mathcal{N}(S^*)^\perp}](y)=P_{\mathcal{N}(S^*)}(y) \in \mathcal{N}(S^*). \] Thus $y \in \cD$ and hence $C^\ominus \subset \cD$. 

For the reverse inclusion, suppose that $y \in \cD$. Then $y=z_0 - (S^*)^\dagger(\rho)$ for some $z_0 \in \mathcal {N}(S^*)$ and $\rho \in \mathcal{N}(S)^\perp$ with $\rho \ge 0$. Then, for each $j=1, 2, \dots, m$, we have \begin{eqnarray*} \la y, k_j \ra &=& \la z_0, k_j \ra - \la (S^*)^\dagger(\rho), k_j\ra=\la z_0, S(e_j) \ra- \la (S^*)^\dagger(\rho), k_j \ra \\ &=& \la S^*(z_0), e_j \ra - \la (S^*)^\dagger(\rho), k_j \ra =- \la (S^*)^\dagger(\rho), S(e_j)\ra \\ &=&-\la S^*(S^*)^\dagger(\rho), e_j\ra = -\la P_{\mathcal{N(S)^\perp}}(\rho), e_j \ra \mbox{  (using Fact \ref{fact 1}(3))}\\ &=& -\la \rho, e_j\ra= -\rho_j \le 0, \end{eqnarray*} which implies that $y\in C^\ominus$ and hence $\cD \subset C^\ominus$.  Thus $\cD=C^\ominus$ and the proof is complete.   \qed 

\medskip Based on this lemma, we can now give a detailed description of best approximations from $C$ and $C^\ominus$ to any $x\in H$. 

\begin{thm}\label{zig char} Let $C$ and $S$ be defined as in equations (\ref{eq 1: fg cone}) and (\ref{eq 2: fg cone}). For each $x\in H$, let $x_0:=x-(S^*)^\dagger S^*(x)$.  Then $x_0 \in \cN(S^*)$ and there exist $\rho, \; \eta \in \R^m$ such that  \begin{enumerate}  \item[{\rm (1)}] $x=S(\rho) +x_0-(S^*)^\dagger(\eta)$. \item[{\rm (2)}]  $\rho \ge 0, \;  \eta \ge 0$, \, $\eta \in \cN(S)^\perp$,\; and  \; $\la \rho, \eta \ra =0$.   \item[{\rm (3)}]  $P_{C^\ominus}(x)=x_0-(S^*)^\dagger (\eta)$,\; and \;$\la x_0, (S^*)^\dagger (\eta) \ra =0$. \item[{\rm (4)}] $P_C(x)=S(\rho)=(S^*)^\dagger[S^*(x) +\eta]$. \end{enumerate} \end{thm} 

\pf  Using Fact \ref{fact 1}(5), we see that \[ S^*(x_0)=S^*(x)-S^*(S^*)^\dagger S^*(x)=S^*(x)-S^*(x)=0 \] and hence $x_0\in \cN (S^*)$. 

 By Fact \ref{moreau}, we have that $x=P_C(x)+P_{C^\ominus}(x)$ and $\la P_C(x), P_{C^\ominus}(x) \ra =0$.  By definition of $C$, $P_C(x)=S(\rho)$ for some $\rho\in \R^m$ with $\rho \ge 0$.  Also, since $P_{C^\ominus}(x)\in C^\ominus$, we use Lemma \ref{zig} to obtain that $P_{C^\ominus}(x)=y-(S^*)^\dagger (\eta)$ for some $y\in \cN(S^*)$ and $-(S^*)^\dagger(\eta)\in \cN(S^*)^\perp$ for some  $\eta\in \cN(S)^\perp$ with  $\eta \ge 0$, and $\la y,  (S^*)^\dagger(\eta) \ra =0$.  We can rewrite this as \be\label{eq 1: zig char} P_{C^\ominus}(x)=x_0 -(S^*)^\dagger (\eta) +y-x_0. \ee Also observe that  \be\label{eq 2: zig char} x-P_{C^\ominus}(x)=\underbrace{(S^*)^\dagger S^*(x)+(S^*)^\dagger(\eta)}_{\in \cN(S^*)^\perp} +\underbrace{x_0-y}_{\in \cN(S^*)} \ee since $(S^*)^\dagger(\eta) \in \cN(S^*)^\perp$ and $S^*(x) \in \cR(S^*)=\cN(S)^\perp$.   \textbf{Claim:}  $y=x_0$. 

For if $y\ne x_0$, then by the Pythagorean theorem we obtain \begin{eqnarray*} \|x-P_{C^\ominus}(x)\|^2&=& \|(S^*)^\dagger S^*(x) +(S^*)^\dagger (\eta) \|^2 +\| x_0-y\|^2\\ &>& \|(S^*)^\dagger S^*(x) +(S^*)^\dagger (\eta) \|^2= \|x-z\|^2, \end{eqnarray*} where $z:=x_0-(S^*)^\dagger(\eta) \in C^\ominus$. This shows that $z$ is a better approximation to $x_0$ from $C^\ominus$ than  $P_{C^\ominus}(x)$ is, which is absurd and proves the claim. 

Thus $P_{C^\ominus}(x)=x_0 -(S^*)^\dagger (\eta)$ and this proves statement (3). Altogether we have that $x=S(\rho)+x_0-(S^*)^\dagger (\eta)$ and this proves statement (1). Statement (4) follows from (3) and Fact \ref{moreau}:  $P_C(x)=x-P_{C^\ominus}(x)$. To verify statement (2), it remains to show that $\la \rho, \eta \ra=0$. But \begin{eqnarray*} 0&=&\la P_C(x), P_{C^\ominus}(x) \ra =\la S(\rho), x_0-(S^*)^\dagger (\eta) \ra \\ &=& \la \rho, S^*(x_0)-S^*(S^*)^\dagger(\eta)\ra \\ &=& \la \rho, -P_{\cN(S)^\perp}(\eta) \ra  \mbox{ \quad (using $x_0\in \cN(S^*)$ and Fact \ref{fact 1}(3))}\\ &=& \la \rho, -\eta \ra  \mbox{ \quad (since $\eta \in \cN(S)^\perp$). } \end{eqnarray*} This completes the proof.   \qed 

\begin{rem} Related to the work of this section, we should mention that Ek\'art, N\'emeth, and N\'emeth \cite{enn;10} have suggested a ``heuristic''  algorithm for computing best approximations from finitely generated cones, in the case where the generators are linearly independent. While they did not prove the convergence of their algorithm, they stated that they numerically solved an extensive set of examples which seemed to suggest that their algorithm was both fast and accurate. 

We believe that Theorems \ref{deu thm}, \ref{deu poly},  and \ref{zig char} will assist us in obtaining an \emph{efficient} algorithm for the actual computation of best approximations from finitely generated cones in Hilbert space. This will be the subject of a future paper.  \end{rem} 

\subsection{An Application to Shape-Preserving Approximation}\label{S: appl} 

In this section, we give a class of problems related to ``shape-preserving''  approximation that can be handled by Theorem \ref{deu thm}. 

Given $x \in L_2[-1,1]$, we want to  find its best approximation from the set of polynomials of degree at most $n$  whose $r$th derivative in nonnegative: \be\label{appl: eq 1} C=C_{n, r}:=\{ p \in \cP_n \mid p^{(r)}(t) \ge 0    \mbox{  for all $t \in [-1,1]$} \}. \ee It is not hard to show that $C$ is a closed convex cone in $L_2[-1,1]$. The interest in such a set is to preserve certain \emph{shape} features of the function being approximated. For example, if $r=0, 1, $ or $2$, then $C$ represents all polynomials of degree $\le n$ that are nonnegative, increasing, or convex, respectively, on $[-1,1]$. It is natural, for example, to want to approximate a convex function in $L_2[-1,1]$ by a \emph{convex} polynomial in $\cP_n$     Choose an orthonormal basis $\{p_0, p_1, \dots, p_n\}$ for $\cP_n$. For definiteness, suppose these are the (normalized) Legendre polynomials. The first five Legendre polynomials are given by  \begin{enumerate} \item[{\rm (0)}] $p_0(t)=\frac{\sqrt{2}}2$,  \item[{\rm (1)}] $p_1(t)=\frac{\sqrt{6}}{2} t$, \item[{\rm (2)}] $p_2(t)=\frac{\sqrt{10}}{4}(3t^2-1)$, \item[{\rm (3)}] $p_3(t)=\frac{\sqrt{14}}{4}(5t^3-3t)$, \item[{\rm (4)}] $p_4(t)=\frac{3\sqrt{2}}{16}(35t^4-30t^2+3)$. \end{enumerate} Thus for each $p\in \cP_n$ we can write its Fourier expansion as $p=\sum_0^n \la p, p_i\ra p_i$.  

For each $\alp \in [-1,1]$, define  \be\label{appl: eq 2} k_\alp:=\sum_{i=0}^n p_i^{(r)}(\alp)p_i \ee  and set   \be\label{Ft: eq 1}  K:=\{k_\alp \mid \alp \in [-1,1]\}.   \ee     \begin{lem}\label{reps} For each $\alp\in [-1,1]$ and $p\in \cP_n$, we have  \be\label{reps: eq 1}  \la k_\alp, p\ra= p^{(r)}(\alp).  \ee  In other words, $k_\alp$ is the representer of the linear functional ``the $r$th derivative evaluated at $\alp$'' on the space $\cP_n$.  \end{lem}    \pf Using the orthonormality of the $p_i$, we get  \begin{eqnarray*}  \la k_\alp, p\ra &=& \left\la \sum_{i=0}^n p_i^{(r)}(\alp)p_i, \sum_{j=0}^n \la p, p_j\ra p_j \right\ra        = \sum_{i=0}^n \sum_{j=0}^n p_i^{(r)}(\alp)\la p, p_j \ra \la p_i, p_j\ra \\       &=&\sum_{i=0}^n \la p, p_i\ra p_i^{(r)}(\alp) =p^{(r)}(\alp).  \hspace{2.5in}\blacksquare \\  \end{eqnarray*} 

 \begin{lem}\label{alt C} \begin{enumerate}  \item[{\rm(1)}] $K$ is a compact set in $\cP_n$.   \item[{\rm(2)}] If $e(t)=t^r$, then $\la e, k\ra =r! > 0$ for all $k\in K$. \item[{\rm(3)}] If $C=C_{n, r}$ is defined as in eq. {\rm (\ref{appl: eq 1})}, then  \be\label{alt C: eq 1}  C=\{ p \in \cP_n \mid \la p, k_\alp \ra \ge 0  \mbox{   for all $\alp\in [-1,1]$} \}.  \ee  \end{enumerate}  \end{lem}    \pf (1) Let $(x_m)$ be a sequence in $K$.  Then there exist $\alp_m \in [-1,1]$ such that $x_m=k_{\alp_m}$ for each $m$. Since the $\alp_m$ are bounded, there is a subsequence $\alp_{m'}$ which converges to some point $\alp \in [-1,1]$. Since $k_{\alp}$ is a continuous function of $\alp$, it follows that  $k_{\alp_{m'}}$ converges to $k_{\alp}$. Thus $K$ is compact.    (2) The $r$th derivative of $t^r$ is the constant $r!$.    (3) This is an immediate consequence of Lemma \ref{reps}.  \qed    \medskip  The following result was first proved in the unpublished thesis of the first  author \cite[Theorem 17]{deu;65}.    \begin{thm}\label{charba in Pn} Let $r, n$ be integers with $0\le r < n$, $X=\cP_n$,  and  \be\label{charba in Pn: eq 1}  C=C_{n, r}:=\{ p \in \cP_n \mid p^{(r)}(t) \ge 0    \mbox{  for all $t \in [-1,1]$} \}.  \ee   Let  $x\in X\setminus C$, $x_0\in C$, and let $k_\alp$ be defined as in {\rm(\ref{appl: eq 2})}. Then the following statements are equivalent:  \begin{enumerate}  \item[{\rm (1)}] $x_0=P_C(x)$;  \item[{\rm (2)}] $x_0=x+\sum_1^m \rho_i k_{\alp_i}$, where $m\le n+1$, $\rho_i > 0$, $\alp_i \in [-1,1]$ and $x_0^{(r)}(\alp_i)=0$ for all $i$, and $\{k_{\alp_1}, k_{\alp_2}, \dots, k_{\alp_m}\}$ is linearly independent.  \end{enumerate}    Moreover, if $x_0^{(r)} \not\equiv 0$ in any of the statements above, then  \[ m\le \frac12(n-r+2). \]   \end{thm}    \pf The equivalence of the statements (1) and (2) is a consequence of Theorem \ref{deu thm} along with  Lemmas \ref{reps}  and  \ref{alt C}.      It remains to show that $m \le (1/2)(n-r+2)$ when $x_0^{(r)}\not\equiv 0$.   Since the vectors $k_{\alp_i}$ are linearly independent, it follows that $\alp_1 \dots, \alp_m$ are distinct points in $[-1, 1]$. Now $x_0^{(r)}$ is a (nonzero) polynomial of degree at most $n-r$, so it has at most $n-r$ zeros. Since $x_0^{(r)}(\alp_i)=0$ for $i=1, \dots, m$, we must have $m\le n-r$. If $x_0^{(r)}(\alp)=0$ for some $\alp$ with $-1< \alp < 1$, then $\alp$ cannot be a simple zero of  $x_0^{(r)}$ (i.e., $x_0^{(r+1)}(\alp)=0$ also) since $x_0^{(r)}(t) \ge 0$ for all $-1\le t \le 1$. It follows that $x_0^{(r)}$ can have at most $1/2(n-r)$ zeros in the open interval $(-1, 1)$. If $x_0^{(r)}$ has a zero at one of the end points $t=\pm 1$, then $x_0^{(r)}$ can have at most $1+1/2(n-r-1)=1/2(n-r+1)$ zeros in $[-1,1]$. Finally, if $x_0^{(r)}$ has zeros at both end points $t=\pm 1$, then we see that $x_0^{(r)}$ has at most $2+1/2(n-r-2)=1/2(n-r+2)$ zeros in $[-1,1]$.  In all possible  cases, we see th  at $x_0^{(r)}$ has at most $m \le 1/2(n-r+2)$ zeros in $[-1,1]$.  \qed             

\section{Elements Vanishing Relative to a Set}\label{S: van char} 

 \begin{defn}\label{vanish} Let $X$ be a normed linear space and $\Gamma \subset X^*$. An element $x^* \in X^*$ is said to {\bf vanish relative to $\Gamma$} if $x\in X$ and $y^*(x)=0$ for all $y^* \in \Gamma$ imply that $x^*(x)=0$.  \end{defn}    Again, when $X$ is a Hilbert space, the above definition reduces to the following form.      \begin{defn}\label{H vanish} Let $X$ be a Hilbert space and $\Gamma \subset X$.  An element $x \in X$ is said to \textbf{vanish relative to $\Gamma$} if $z\in X$ and $\la z, y\ra =0$ for all $y\in \Gamma$ imply that $\la x, z \ra =0$.  \end{defn}    This idea can be characterized in a useful way just as ``positive relative to a set'' was in Theorem \ref{pos char}.     \begin{thm}\label{van char} Let $X$ be a normed linear space, $\Gamma \subset X^*$, and $x^* \in X^*$.  Then the following statements are equivalent:  \begin{enumerate}  \item[{\rm (1)}]  $x^*$ vanishes relative to $\Gamma$.   \item[{\rm (2)}] $\Gamma_\perp \subset (x^*)_\perp$.  \item[{\rm (3)}] $x^* \in w^*\!\!-cl(\spane(\Gamma))$, the weak* closed linear span of $\Gamma$.    Moreover, if $X$ is reflexive, then each of these statements is equivalent to 

 \indent \!\!\!\!\!\!\!\!\!{\rm (4)} $x^* \in \overline{\spane}(\Gamma)$, the (norm) closed linear span of $\Gamma$.  \end{enumerate}  \end{thm}    \pf $(1) \Rightarrow (2)$. Suppose (1) holds. Let $x \in \Gamma_\perp$. Then $y^*(x)=0$ for all $y^* \in \Gamma$.  By (1), $x^*(x)=0$. That is, $x \in (x^*)_\perp$.  Hence (2) holds.    $(2) \Rightarrow (3)$. If (3) fails, $x^* \notin w^*\!\!-cl(\spane(\Gamma))$. By Theorem \ref{big thm}, there exists a weak$^*$ continuous linear functional $f$ on $X^*$ such that   \be\label{van char: eq 1}  \sup f(\spane(\Gamma)) < f(x^*).  \ee  But $f=\hat{x}$ for some $x\in X$. Thus we can rewrite the inequality (\ref{van char: eq 1}) as  $\sup \hat{x}(\spane(\Gamma)) < \hat{x}(x^*)$, or  \be\label{van char: eq 2}  \sup\{ y^*(x) \mid y^* \in \spane(\Gamma) \} < x^*(x).  \ee  Since $\spane(\Gamma)$ is a linear subspace, the only way the expression on the left side of (\ref{van char: eq 2}) can be bounded above is if $y^*(x)=0$ for each $y^* \in \Gamma$. In this case, it follows that $x^*(x)>0$.  Thus $x\in \Gamma_\perp \setminus (x^*)_\perp$ and (2) fails.    $(3) \Rightarrow (1)$. Let $x^* \in w^*\!\!-cl(\spane(\Gamma))$.  If $x\in X$ and $y^*(x)=0$ for all $y^*\in \Gamma$, then clearly $y^*(x)=0$ for all $y^* \in \spane(\Gamma)$. Since $x^* \in w^*\!\!-cl(\spane(\Gamma))$, there exists a net $(y_\alp^*) \in \spane(\Gamma)$ such that $y^*_\alp$ weak$^*$ converges to $x^*$, i.e., $y^*_\alp(z) \to x^*(z)$ for each $z\in X$. But $y_\alp^*(x)=0$ for all $\alp$, so $x^*(x)=0$. That is, $x^*$ vanishes relative to $\Gamma$ and (1) holds.       If $X$ is reflexive, then the same proof as in Theorem \ref{pos char} works.  \qed     \begin{cor}\label{pos implies van} Let $X$ be a normed linear space, $\Gamma \subset X^*$, and $x^*\in X^*$. If $x^*$ is positive relative to $\Gamma$, then $x^*$ vanishes relative to $\Gamma$.  \end{cor}    \pf  By Theorem \ref{pos char}, we have $x^* \in w^*\!\!-cl(\cone(\Gamma))$. A fortiori, $x^*\in w^*\!\!-cl (\spane(\Gamma))$. By Theorem \ref{van char}, the result follows.     A simpler, more direct, proof goes as follows. For any subset $S$ of $X^*$, we have $S_\ominus \cap (-S_\ominus)=S_\perp$.  Hence if $\Gamma_\ominus \subset (x^*_\ominus)$, then $-\Gamma_\ominus \subset -(x^*_\ominus)$ and hence  $$  \Gamma_\perp=\Gamma_\ominus \cap[-\Gamma_\ominus]\subset (x^*)_\ominus \cap [-(x^*)_\ominus]=(x^*)_\perp  $$  In other words, using Theorem \ref{van char}, $x^*$ vanishes relative to $\Gamma$.  \qed    \medskip The following simple example shows that the converse to this theorem is not valid.    \begin{ex}\label{van not pos} Let $X=\R$ with the absolute-value norm: $\|x\|:=|x|$. Let $x=-1$ and $\Gamma=\{1\}$. Then $x$ vanishes relative to $\Gamma$, but $x$ is not positive relative to $\Gamma$.   \end{ex}    \medskip     Again, by the same argument as in Theorem \ref{pos char}, we note that in a reflexive Banach space $X$, a convex set in $X^*$ is (norm) closed if and only if it is weak$^*$ closed.  Thus we have the following result.    \begin{cor}\label{R cor van} Let $X$ be a reflexive Banach space, $\Gamma \subset X^*$, and $x^* \in X^*$. Then $x^*$ vanishes relative to $\Gamma$ if and only if $x^* \in \overline{\spane}(\Gamma)$, the (norm) closed linear span of $\Gamma$.  \end{cor} 

\begin{cor}\label{H cor van} Let $H$ be a Hilbert space, $\Gamma \subset H$, and $x \in H$.  Then $x$ vanishes relative to $\Gamma$ if and only if $x \in \overline{\spane}({\Gamma})$, the norm closed linear span of $\Gamma$. \end{cor}

One important application of Theorem \ref{van char} is the following. 

\begin{lem}\label{kernels} Let $X$ be a normed linear space and $f, f_1, \dots, f_n$ be in $X^*$.  Then the following statements are equivalent: \begin{enumerate} \item[{\rm (1)}] $f$ vanishes relative to $\Gamma=\{f_1, f_2,\dots, f_n\}$. \item[{\rm (2)}] If $x \in X$ and $f_i(x)=0$ for each $i=1,2, \dots, n$, then $f(x)=0$ . \item[{\rm (3)}]  $\cap_{i=1}^n f_i^{-1}(0) \subset f^{-1}(0)$. \item[{\rm (4)}] $f\in \spane\{f_1, f_2, \dots, f_n\}$. \end{enumerate} \end{lem} 

\pf  The equivalence of (1) and (2) is just a rewording of the definition, and the equivalence of (2) and (3) is obvious.  Finally, (1) holds if and only if $f$ lies in the weak$^*$ closed linear span of $\Gamma$ by Theorem \ref{van char}.  But the linear span of $\Gamma$, being finite-dimensional, is weak$^*$ closed (see, eg., \cite[Corollary 3.14]{fhhmz;11}). That is, (4) holds.  \qed 

\medskip  This result---in even more general vector spaces---has proven useful in studying weak topologies on vector spaces (see, e.g., \cite[p. 421]{dusc;58} or \cite[Lemma 3.21]{fhhmz;11}).  

\section*{Acknowledgements} The work of L.~Zikatanov was supported in part by NSF awards DMS-1720114 and DMS-1522615. 

                                     \bibliographystyle{plain}  

\begin{tabular}{lll} Frank Deutsch    &  Hein Hundal \\ Department of Mathematics  &   146 Cedar Ridge Drive  \\ Penn State University  & Port Matilda, PA 16870 \\ University Park, PA 16802  & email: hundalhh@yahoo.com \\  email: deutsch@math.psu.edu & { } \\  {} &{}\\   Ludmil Zikatanov & {} \\   Department of Mathematics &{}\\  Penn State University &{}\\  University Park, PA 16802 &{}\\  email:  ltz@math.psu.edu &{}\\ \end{tabular} 

\end{document}